\documentclass[11pt]{amsart}
\theoremstyle{plain}
\newtheorem{thm}{Theorem}[section]
\newtheorem{theorem}[thm]{Theorem}

\newtheorem{lemma}[thm]{Lemma}
\newtheorem{corollary}[thm]{Corollary}
\newtheorem{proposition}[thm]{Proposition}
\theoremstyle{definition}

\newtheorem{definition}[thm]{Definition}

\newtheorem{example}[thm]{Example}

\newtheorem{conjecture}[thm]{Conjecture}

\newtheorem{question}[thm]{Question}

\numberwithin{equation}{section}

\newcommand{\p}{\partial}

\newcommand{\sC}{{\mathcal C}}

\newcommand{\sF}{{\mathcal F}}

\newcommand{\sH}{{\mathcal H}}

\newcommand{\sJ}{{\mathcal J}}
\newcommand{\sK}{{\mathcal K}}

\newcommand{\sO}{{\mathcal O}}
\newcommand{\sP}{{\mathcal P}}
\newcommand{\sQ}{{\mathcal Q}}

\newcommand{\sU}{{\mathcal U}}

\newcommand{\sZ}{{\mathcal Z}}

\newcommand{\C}{{\mathbb C}}

\newcommand{\F}{{\mathbb F}}

\newcommand{\BP}{{\mathbb P}}

\newcommand{\Q}{{\mathbb Q}}

\newcommand{\BS}{{\mathbb S}}

\newcommand{\w}{\widetilde}
\newcommand{\h}{\widehat}

\newcommand{\fg}{{\mathfrak g}}

\newcommand{\fgl}{{\mathfrak g}{\mathfrak l}}

\def\Gr{\mathop{\rm Gr}\nolimits}

\def\Sym{\mathop{\rm Sym}\nolimits}

\def\Hom{\mathop{\rm Hom}\nolimits}

\title[Neighborhoods of Minimal Rational Curves]{Geometry of Neighborhoods of Minimal Rational Curves}

\author{Jun-Muk Hwang}

\thanks{This work was supported by the Institute for Basic Science (IBS-R032-D1).}
\begin{document}

\begin{abstract}
This is a survey of recent works on the germ-equivalence problem of minimal rational curves on uniruled projective manifolds. 
Our main interest is  when the associated varieties of minimal rational tangents form an isotrivial family of projective varieties. In this case, there is  a natural G-structure on a Zariski-open subset of the underlying uniruled projective manifold, which leads to an interaction of algebraic geometry of minimal rational curves with differential geometry of geometric structures.   We also discuss the related question of the formal principle for the germ-equivalence of minimal rational curves. 
\end{abstract}

\maketitle




\section{Introduction: germ-equivalence of minimal rational curves}\label{germ} 
We work in the setting of complex analytic geometry. All manifolds, varieties and morphisms are complex-analytic. Open sets are those in the Euclidean topology. Open sets in the Zariski topology are called Zariski-open sets. A nonsingular projective variety is called a projective manifold. The projectivization $\BP V$ of a vector space  $V$ is the set of one-dimensional subspaces in $V$. 

\medskip
Let $X$ be a uniruled projective manifold.  For an irreducible component $\sK$ of the space of rational curves on $X$, let 
  \begin{equation}\label{eq.mu} \sK \ \stackrel{\rho}{\longleftarrow} \  {\rm Univ}_{\sK} \ \stackrel{\mu}{\longrightarrow} \ X \end{equation} be its universal family: for a point $z \in \sK$, the fiber $\rho^{-1}(z)$ is isomorphic to $\BP^1$ and $\mu|_{\rho^{-1}(z)}$ is the normalization of the rational curve $\mu(\rho^{-1}(z)) \subset X$ (see \cite[Section II.2]{Kol}). We call $\sK$  a {\em family of minimal rational curves}, if the fiber $\mu^{-1}(x)$ for a general point $x \in X$ is nonempty and projective.  A member of $\sK$ is called a {\em minimal rational curve } on $X$. For example, when we have a projective embedding $X \subset \BP^N,$ rational curves of minimal degree through general points of $X$ are minimal rational curves.  We are interested in the following equivalence of germs of minimal rational curves. 

\begin{definition}\label{d.germ}
Let $X$ and $\w{X}$ be uniruled projective manifolds. Minimal rational curves $C \subset X$ and $\w{C} \subset \w{X}$ are {\em germ-equivalent}, if there exist  neighborhoods $C \subset U \subset X$ and $\w{C} \subset \w{U} \subset \w{X}$ equipped with a biholomorphic map $\varphi: U \to \w{U}$ satisfying $\varphi(C) = \w{C}$. \end{definition}

One motivation to study the above notion is the following result from \cite{HM01}, a version of Cartan-Fubini type extension theorem.

\begin{theorem}\label{t.CF}
Let  $X$ and $\w{X}$ be uniruled projective manifolds of Picard number 1. Let $\sK$ (resp. $\w{\sK}$) be a family of minimal rational curves on $X$ (resp. $\w{X}$) such that  the  universal family morphism $\mu: {\rm Univ}_{\sK} \to X$ (resp. $\w{\mu}: {\rm Univ}_{\w{\sK}} \to \w{X}$) has connected fibers. If  a member $C \subset X $ of $\sK$ is germ-equivalent to a member $\w{C} \subset \w{X}$ of $\w{\sK}$ by a biholomorphic map $\varphi:U \to \w{U}$ between some neighborhoods,  then $\varphi$ can be extended to a biregular morphism from  $X$ to $\w{X}$. \end{theorem} 

Theorem \ref{t.CF} has been useful in proving certain rigidity results on a class of Fano manifolds of Picard number 1. This shows that it is worthwhile to study  the germ-equivalence of minimal rational curves. 
Note that if $C \subset X$ is germ-equivalent to $\w{C} \subset \w{X}$ via a biholomorphic map $\varphi: U \to \w{U}$, then a deformation of $C$ inside the neighborhood $U$ is germ-equivalent to some  deformation of $\w{C}$ inside the neighborhood $\w{U}$. 
Thus when our interest is on applications like Theorem \ref{t.CF}, it is sufficient to consider the case when $C \subset X$ (resp. $\w{C} \subset \w{X}$) is a general member of $\sK$ (resp. $\w{\sK}$). By Mori's bend-and-break argument, a general member 
 $C \subset X$ of $\sK$ is  an {\em unbendable rational curve}, namely, a  rational curve whose normalization $\nu: \BP^1 \to C \subset X$ satisfies \begin{equation}\label{e.N} \nu^* TX  \ \cong \ \sO(2) \oplus \sO(1)^d \oplus \sO^{n-1-d}, \ n := \dim X,\end{equation} where $d$ is the dimension of the fiber of the universal family morphism $\mu: \sK \to X$.  Furthermore, when $d >1$, a general member of $\sK$ is  a smooth rational curve (\cite[Theorem 3.3]{Ke}). 
Thus when we study the germ-equivalence of minimal rational curves, it is reasonable to concentrate on the case of  smooth unbendable minimal rational curves. Obviously, two germ-equivalent unbendable rational curves have the same values of $n$ and $d$ in  (\ref{e.N}). When $d=0$, the normal bundle of the rational curve is trivial. It is easy to see that two smooth rational curves with trivial normal bundles of the same rank are always germ-equivalent. So we may assume that $d>0$.   Then the basic question for us is the following.

\begin{question}\label{q.equiv}
Fix positive integers $n$ and $d <n$. How large is the set of germ-equivalence classes of  smooth minimal rational curves $C \subset X$  satisfying (\ref{e.N})? Can we classify them? \end{question}

I claim that the set of germ-equivalence classes in Question \ref{q.equiv} is huge, bigger than the set of the following equivalence classes in projective algebraic geometry.

\begin{definition}\label{d.proj}
Let $Z, \w{Z} \subset \BP^{n-1}$  be two $d$-dimensional  projective submanifolds of the complex projective space of dimension $n-1$. We say that $Z$ and $\w{Z}$ are {\em projectively equivalent} if there exists a projective linear automorphism $\psi \in {\rm PGL}(\C^n)$ of $\BP^{n-1}$ such that $\psi(Z) = \w{Z}$. \end{definition} 

What is the relation between the equivalence classes in Definition \ref{d.germ} and those in Definition \ref{d.proj}?  It is easy to see that  we can associate to each projective submanifold $Z \subset \BP^{n-1}$ a family of minimal rational curves  in the following way. 

\begin{proposition}\label{p.blowup}
Given a projective submanifold $Z \subset \BP^{n-1}$ of dimension $d$,   regard $\BP^{n-1}$ as a hyperplane in $\BP^n$ and let $\beta: X \to \BP^n$  be the blowup of $\BP^n$ along the submanifold $Z \subset \BP^{n-1} \subset \BP^n$.
Let $\sK^o$  be the set  of rational curves $C \subset X$  satisfying \begin{itemize} \item[(1)] $\beta(C)$ is a line on $\BP^{n}$; \item[(2)] $\beta|_C: C \to \beta(C)$ is biregular; \item[(3)] $\beta(C) \not\subset \BP^{n-1}$; and
\item[(4)] $\beta(C) \cap Z \neq \emptyset.$ \end{itemize}
Then the closure $\sK$ of $\sK^o$ in the space of rational curves on $X$ is a family of minimal rational curves on $X$ and a member $C$ of $\sK^o$  is a    smooth unbendable minimal rational curve satisfying (\ref{e.N}). \end{proposition}
 
The following proposition says that the set of equivalence classes in Question \ref{q.equiv} is larger than that of Definition \ref{d.proj}. 

\begin{proposition}\label{p.huge}
Given two projective submanifolds $Z, \w{Z} \subset \BP^{n-1}$, let $\sK^o \subset \sK$ (resp. $\w{\sK}^o \subset \w{\sK}$) be a family of minimal rational curves on a uniruled projective manifold $X$ (resp. $\w{X}$) obtained by applying Proposition \ref{p.blowup} to $Z$ (resp. $\w{Z}$) and let $C \subset X$ (resp. $\w{C} \subset \w{X}$) be a member of $\sK^o$ (resp. $\w{\sK}^o$).   If $C \subset X$ and $\w{C} \subset \w{X}$ are germ-equivalent, then $Z$ and $\w{Z}$ are projectively equivalent.  \end{proposition}

To see Proposition \ref{p.huge},  we need the following notion of VMRT (an abbreviation of Variety of Minimal Rational Tangents). We recommend \cite{HM99}, \cite{Hw01} and \cite{Hw11} for introductory surveys on VMRT. 

\begin{definition}\label{d.vmrt}
Let $\sK$ be a family of minimal rational curves on a projective manifold $X$. Sending a nonsingular point $x \in C$ on a member of $\sK$ to its tangent direction $\BP T_x C \in \BP T_x X$ defines a rational map $\tau: {\rm Univ}_{\sK} \dasharrow \BP TX$ to the projectivized tangent bundle of $X$.  The closure of the image of $\tau$ is an irreducible projective subvariety $\sC \subset \BP TX$, called the {\em total VMRT} of $\sK$, and the map $\tau: {\rm Univ}_{\sK} \dasharrow \sC$ is birational by \cite[Theorem 1]{HM04}. For a general point $x \in X$, the fiber of $\sC \to X$ is a $d$-dimensional projective subvariety $\sC_x \subset \BP T_x X$ (with finitely many components), called the {\em VMRT } of $\sK$ at $x$. \end{definition}

   Let $C \subset X$ be a smooth unbendable minimal rational curve. By (\ref{e.N}), for each point $x \in C$, there is a $d$-dimensional family of deformations of $C$ fixing $x$, which are also unbendable rational curves. Their tangent vectors at $x$ determine a $d$-dimensional locally closed submanifold in $\BP T_x M$ which lies in the VMRT $\sC_x$.   
From this, the following lemma is immediate. 

\begin{lemma}\label{l.vmrt}
In Definitions \ref{d.germ} and \ref{d.vmrt}, assume that the VMRT at a general of $X$ (and $\w{X}$)  is smooth and irreducible. If $C \subset X$ and $\w{C} \subset \w{X}$ are germ-equivalent by 
  a biholomorphic map $\varphi: U \to \w{U}$ between some neighborhoods, then for each $x \in C$ and $\w{x} := \varphi(x) \in \w{C}$, 
the two projective submanifolds $$\sC_x \subset \BP T_x M \cong \BP^{n-1} \mbox{ and } \w{\sC}_{\w{x}} \subset \BP T_{\w{x}} \w{M} \cong \BP^{n-1}$$ are projectively equivalent by the projective linear isomorphism $${\rm d}_x \varphi: \BP T_x M \to \BP T_{\w{x}} \w{M}.$$ \end{lemma} 

Then the proof of Proposition \ref{p.huge}  follows from the next lemma, which is easy to check. 

\begin{lemma}\label{l.blowup}
In Proposition \ref{p.blowup},  for any point $x \in X \setminus \beta^{-1}(\BP^{n-1})$, the VMRT $ \sC_x  \subset \BP T_x X$ at $x$ of the family $\sK$ is projectively equivalent to $Z \subset \BP^{n-1}$. \end{lemma}

Proposition \ref{p.huge} shows that the classification in Question \ref{q.equiv} requires at least a classification of projective equivalence classes in Definition \ref{d.proj} of all projective submanifolds of dimension $d$ in $\BP^{n-1}$. The latter is, in some sense, the ultimate goal of  projective algebraic geometry! Thus Question \ref{q.equiv} as it stands is too big a problem. To make it a reasonable problem, we have to take into account the vast possibilities of germ-equivalence classes of minimal rational curves arising from Proposition \ref{p.blowup}. In this regard, it is convenient to introduce the following terminology.

\begin{definition}\label{d.typeS}
A minimal rational curve $C \subset X$ with (\ref{e.N}) is said to be {\em of locally flat type }  if it is germ-equivalent to some $C \subset X$ in Proposition \ref{p.blowup} for a suitable choice of a $d$-dimensional projective submanifold $Z \subset \BP^{n-1}$.\end{definition}

From Lemma \ref{l.blowup}, a necessary condition for a minimal rational curve $C \subset X$ to be of locally flat type  is that the VMRT $ \sC_x \subset \BP T_x X$ at a general point $x \in X$ must be projectively equivalent to a fixed nonsingular projective variety $Z \subset \BP^{n-1}$. It is natural to ask the converse:

\begin{question}\label{q.Z}
Fix a positive-dimensional projective submanifold $Z \subset \BP^{n-1}$.
Let $\sK$ be a family of minimal rational curves on a projective manifold $X$ such that   \begin{center} $(\dagger)$  the VMRT 
 $\sC_x \subset \BP T_x X$ at a general point $x \in X$ is projectively equivalent to $Z \subset \BP^{n-1}$. \end{center} Is a general member $C \subset X$ of $\sK$ of locally flat type? If this is not the case, can we classify the germ-equivalence classes of smooth unbendable minimal rational curves satisfying the assumption $(\dagger)$? \end{question}
 
 Question \ref{q.Z} is exactly 
 a more reasonable version of Question \ref{q.equiv}, taking into account those examples arising from Proposition \ref{p.blowup}. Question \ref{q.Z} is 
 the main topic of our discussion. It is actually an infinite collection of questions, one independent question for each choice of $Z  \subset \BP^{n-1}$ up to projective equivalence. Our approach to Question \ref{q.Z} is to combine algebraic geometry of minimal rational curves with differential geometry of geometric structures. In particular, the term `locally flat' in Definition \ref{d.typeS} comes from differential geometry, which we explain in the next section.  
 
 \section{$Z$-cone structures}\label{s.cone}
 Here we explain a differential-geometric structure arising from Question \ref{q.Z}.  
  Fix   a vector space $V$ of dimension $n$ and a  projective submanifold $Z \subset \BP V.$ 
 
 \begin{definition}\label{d.cone}
 Let $M$ be a complex manifold of dimension $n$. A closed submanifold $\sC \subset \BP TM$ is a {\em $Z$-cone structure} on $M$ if the fiber $\sC_x \subset \BP T_x M$ for each $x \in M$ of the natural projection $\sC \to M$ is projectively equivalent to $Z \subset \BP V$. \end{definition}

 There is a more general notion of cone structures (see \cite[Section 3]{Hw11} or \cite[Section 1.2]{HN}). The more restricted notion in Definition \ref{d.cone} is sufficient for our purpose here. We mention that a systematic study of cone structures had originated from twistor theory (see \cite[Chapter 1, Section 6]{Ma}). 
 From the general perspective of differential-geometric structures, the following is the natural equivalence notion of cone structures.

 \begin{definition}\label{d.equiv}
 Let $\sC \subset \BP TM$ (resp. $\w{\sC} \subset \BP T \w{M}$) be a $Z$-cone structure on a complex manifold $M$ (resp. $\w{M}$). Then $\sC$ at $x \in M$ is   {\em locally equivalent} to $\w{\sC}$ at $\w{x} \in \w{M}$, if there exist  neighborhoods $x \in O \subset M$ and $\w{x} \in \w{O} \subset \w{M}$ equipped with a biholomorphic map $\phi: O \to \w{O}$ such that $\w{x} = \phi(x)$ and the differential ${\rm d} \phi: \BP T O \to \BP T \w{O}$ sends $\sC|_O$ to $\w{\sC}|_{\w{O}}$, yielding the following commutative diagram:
 $$\begin{array}{ccccccc}
 \sC & \supset & \sC|_{O} & \stackrel{{\rm d} \phi}{\longrightarrow} & \w{\sC}|_{\w{O}} & \subset & \w{\sC} \\ \cap & & \cap & & \cap & & \cap \\
 \BP TM & \supset & \BP T O & \stackrel{{\rm d} \phi}{\longrightarrow} & \BP T \w{O} & \subset & \BP T \w{M} \\ \downarrow & & \downarrow & & \downarrow & & \downarrow \\
 M & \supset & O & \stackrel{\phi}{\longrightarrow} & \w{O} & \subset & \w{M}. \end{array} $$ 
 \end{definition}
 
 Obviously, in Question \ref{q.Z}, there is a Zariski-open subset $M \subset X$ equipped with a $Z$-cone structure $\sC \subset \BP TM$ coming from the VMRT. Definition \ref{d.equiv} is relevant in Question \ref{q.Z} because of the following result from \cite[Proposition 2.1]{HM01}.
 
 \begin{proposition}\label{p.CF}
Fix a  projective submanifold $Z \subset \BP V$, which is not a linear subspace. 
Let $\sK$ (resp. $\w{\sK}$) be a family of minimal rational curves on a projective manifold $X$ (resp. $\w{X}$) such that   the VMRT 
 $\sC_x \subset \BP T_x X$ (resp. $\w{\sC}_{\w{x}} \subset \BP T_{\w{x}} \w{X}$)  at a general point $x \in X$ (resp. $\w{x} \in \w{X}$) is projectively equivalent to $Z \subset \BP V$, giving rise to a $Z$-cone structure $\sC \subset \BP TM$ (resp. $\w{\sC} \subset \BP T\w{M}$) on a Zariski-open subset $M \subset X$ (resp. $\w{M} \subset \w{X}$). Then the following two conditions are equivalent.  \begin{itemize} \item[(i)] For a general member $C \subset X$ of $\sK$, there exists a member $\w{C} \subset \w{X}$ of $\w{\sK}$ such that $C$ and $\w{C}$ are germ-equivalent (in the sense of Definition \ref{d.germ}). \item[(ii)] The $Z$-cone structure $\sC$ at some points $x \in M$ is locally equivalent to the $Z$-cone structure $\w{\sC}$ at some point $\w{x} \in \w{M}$ (in the sense of Definition \ref{d.equiv}). \end{itemize} \end{proposition} 
 
 Note that in \cite{HM01}, the result was stated under the assumption that the projective variety $Z \subset \BP V$ has nondegenerate Gauss map. It implies Proposition \ref{p.CF} because this assumption is satisfied by any projective submanifold which is not a linear subspace. By Proposition \ref{p.CF}, Question \ref{q.Z} can be treated as a problem of local equivalence of $Z$-cone structures, a question in differential geometry.
  
  A standard example of a $Z$-cone structure  is the following flat structure.

\begin{definition}\label{d.flat}
Regard $V$ as a complex manifold. We have a canonical trivialization of the tangent bundle $TV = V \times V$, which induces canonical a trivialization $\BP TV = (\BP V) \times V$. Fix a projective submanifold $Z \subset \BP V$. \begin{itemize} \item[(i)]  The $Z$-cone structure on $V$ 
$$ \sC^Z:= Z \times V \ \subset \ \BP V \times V = \BP TV $$ is called the {\em flat $Z$-cone structure}. 
\item[(ii)] A $Z$-cone structure $\sC \subset \BP TM$ on a complex manifold $M$ is {\em locally flat} if $\sC$ at some point $x \in M$  is locally equivalent  to the flat $Z$-cone structure $\sC^Z$  at a point on $V$. \end{itemize} 
\end{definition}
 
 It is easy to see that the $Z$-cone structure on the Zariski-open subset $M:= X \setminus \beta^{-1}(\BP^{n-1})$ in Proposition \ref{p.blowup} is locally flat. This explains the terminology `locally flat' in Definition \ref{d.typeS} and Question \ref{q.Z}.

 The following  weaker notion than local flatness is also important.
 
 \begin{definition}\label{d.sym}
 A $Z$-cone structure $\sC \subset \BP TM$ is {\em locally symmetric} if a general $x \in M$ admits a neighborhood $O$ with a biholomorphic map $\sigma: O \to \sigma(O) \subset M$ such that $\sigma(x) =x, {\rm d}_x \sigma = - {\rm Id}_{T_x M}$ and ${\rm d} \sigma: \BP TO \to \BP T\sigma(O)$ sends $\sC|_O$ to $\sC|_{\sigma(O)}$. \end{definition}
 
 Locally flat $Z$-cone structures are locally symmetric, but the converse is not true:  see the examples in Subsection  \ref{s.adjoint}. 
 
 A $Z$-cone structure arising from varieties of minimal rational tangents is equipped with some additional structures. To discuss them, we need the following notion.  
  
 \begin{definition}\label{d.conic}
 Let $\sC \subset \BP TM$ be a $Z$-cone structure with the natural projection $\psi: \sC \to M$. \begin{itemize} \item[(i)] For each $\alpha \in \sC$ and $x= \psi(\alpha),$ consider the differential ${\rm d}_{\alpha} \psi: T_{\alpha} \sC \to T_x M$ and the 1-dimensional subspace $\h{\alpha} \subset T_x M$ corresponding to $\alpha \in \BP T_x M$. Define $$\sJ_{\alpha} \ := \ ({\rm d}_{\alpha} \psi)^{-1} (\h{\alpha}) \ \supset \ {\rm Ker} ({\rm d}_{\alpha} \psi).$$ Then $\sJ := \cup_{\alpha \in \sC} \sJ_{\alpha}$ is a vector subbundle of $T \sC$ with a natural exact sequence $$0 \ \longrightarrow \  {\rm Ker}({\rm d} \psi) \  \longrightarrow \  \sJ \ \longrightarrow \ L \ \longrightarrow 0,$$ where $L$ is the line bundle on $\sC$ defined by the restriction of the tautological line bundle of the projectivization $\BP TM$.  \item[(ii)] A {\em conic connection} on $\sC$ is a line subbundle $\sF \subset \sJ$ which splits the exact sequence in ${\rm (i)}$, giving a direct sum decomposition $\sJ \cong {\rm Ker}({\rm d} \psi) \oplus \sF.$  \item[(iii)] A conic connection $\sF \subset \sJ$ is a {\em characteristic conic connection} if $$[\sF,  [\sJ, \sJ]] \ \subset \  \sJ + [\sJ, \sJ],$$ where we regard $\sF$ and $\sJ$ as locally free sheaves of vector fields on $\sC$, and the brackets denote the Lie brackets of vector fields.   \end{itemize} \end{definition}
  
  A key result is the following from \cite[Proposition 8]{HM04}. 
  
  \begin{proposition}\label{p.conic}
 The $Z$-cone structure on a Zariski-open subset $M \subset X$ arising from Question \ref{q.Z} has a natural conic connection coming from minimal rational curves: the fibers of the morphism  $\rho$ in (\ref{eq.mu}) are sent by the birational map $\tau: {\rm Univ}_{\sK} \dasharrow \sC$ in Definition \ref{d.vmrt} to a line subbundle $\sF \subset \sC$ over $M$. Furthermore, it is a characteristic conic connection. \end{proposition} 
  
   From Propositions \ref{p.CF} and \ref{p.conic},  Question \ref{q.Z} suggests  the following problem in differential geometry. 
  
  \begin{question}\label{q.Zcone}
  Fix a positive-dimensional projective submanifold $Z \subset \BP V$. 
  Let $\sC \subset \BP TM$ be a $Z$-cone structure equipped with a characteristic conic connection.
  Is it locally flat (resp. locally symmetric)? Can we classify the local equivalence classes of such $Z$-cone structures? \end{question}
  
  There are differential-geometric tools to study the local equivalence of $Z$-cone structures, which we discuss in the next section. 
  
  \section{G-structures}\label{s.G}
  
 To study the local equivalence problem for $Z$-cone structures, it is convenient to use
  the theory of G-structures, a very classical subject in differential geometry (see \cite[Chapter I]{Ko} and \cite[Section 2]{Sc}).  Here, we just recall some basic ingredients of the theory, which are necessary for our discussion of $Z$-cone structures. 
  
  \begin{definition}\label{d.pb}
  Let $G \subset H$ be a complex Lie subgroup of a complex Lie group and let $M$ be a complex manifold. \begin{itemize} \item[(i)]  A {\em principal bundle on $M$ with the structure group} $H$ is a complex manifold $\sP$ equipped with 
   a free right $H$-action on $\sP$ and a surjective morphism $\pi: \sP \to M$ whose fibers are $H$-orbits. 
  \item[(ii)] A closed submanifold $\sQ \subset \sP$ is a {\em principal subbundle of $\sP$ with the structure group $G$} if the induced right action of $G \subset H$ on $\sP$ makes $\sQ$ a principal bundle on $M$ with the structure group $G$. \end{itemize} \end{definition}
 
 \begin{definition}\label{d.Gstructure}
 Fix a vector space $V$ of dimension $n$ and a subgroup $G \subset {\rm GL}(V)$. \begin{itemize} \item[(i)]
 For a complex manifold $M$ of dimension $n$, let $\F M := \cup_{x \in M} \F_x M$ be the {\em frame bundle} with the fiber at $x \in M$ given by $$\F_x M \ := \ {\rm Isom}(V, T_x M).$$ It is a principal bundle on $M$ with the structure group ${\rm GL}(V)$.  \item[(ii)] A principal subbundle $\sP \subset \F M$ with the structure group $G$ is called a  {\em $G$-structure on $M$}. \item[(iii)]
  Let $\sP \subset \F M$ (resp. $\w{\sP} \subset \F \w{M}$) be a $G$-structure on a complex manifold $M$ (resp. $\w{M}$). Then $\sP$ at $x \in M$ and $\w{\sP}$ at $\w{x} \in \w{M}$ are {\em locally equivalent} if there exist  neighborhoods $x \in O \subset M$ and $\w{x} \in \w{O} \subset \w{M}$ equipped with a biholomorphic map $\phi: O \to \w{O}$ such that  $\w{x} = \phi(x)$ and the differential ${\rm d} \phi: \F O \to \F \w{O}$ sends $\sP|_O$ to $\w{\sP}|_{\w{O}}$, yielding the following commutative diagram:
 $$\begin{array}{ccccccc}
 \sP & \supset & \sP|_{O} & \stackrel{{\rm d} \phi}{\longrightarrow} & \w{\sP}|_{\w{O}} & \subset & \w{\sP} \\ \cap & & \cap & & \cap & & \cap \\
 \F M & \supset & \F O & \stackrel{{\rm d} \phi}{\longrightarrow} & \F \w{O} & \subset & \F \w{M} \\ \downarrow & & \downarrow & & \downarrow & & \downarrow \\
 M & \supset & O & \stackrel{\phi}{\longrightarrow} & \w{O} & \subset & \w{M}. \end{array} $$ \end{itemize}
 \end{definition}
 
  A standard example of a $G$-structure  is the following flat structure.

\begin{definition}\label{d.Gflat}
Regard $V$ as a complex manifold. We have a canonical trivialization of the tangent bundle $TV = V \times V$, which induces a canonical trivialization $\F V = {\rm GL}(V) \times V$. Fix a subgroup $G \subset {\rm GL}(V)$. \begin{itemize} \item[(i)]  The $G$-structure on $V$
$$ \sP^Z:= G \times V \ \subset \ {\rm GL}(V) \times V = \F V $$ is called the {\em flat $G$-structure}. 
\item[(ii)] A $G$-structure $\sP \subset \F M$ on a complex manifold $M$ is {\em locally flat} if $\sP$ at some point $x \in M$  is locally equivalent  to the flat $G$-structure $\sP^Z$ at a point on $V$. \end{itemize} 
\end{definition}
 
 Analogously to Definition \ref{d.sym}, we have the following.
 
 \begin{definition}\label{d.Gsym}
 A $G$-structure $\sP \subset \F M$ is {\em locally symmetric} if a general $x \in M$ admits a neighborhood $O$ with a biholomorphic map $\sigma: O \to \sigma(O)$ such that $\sigma(x) =x, {\rm d}_x \sigma = - {\rm Id}_{T_x M}$ and ${\rm d} \sigma: \F O \to \F \sigma(O)$ sends $\sP|_O$ to $\sP|_{\sigma(O)}$. \end{definition}
 
 To study whether a $G$-structure is locally flat or locally symmetric, it is useful to introduce connections. 
 
  \begin{definition}\label{d.principal}
  Let $G$ be a complex Lie group. 
A {\em principal connection} on a principal bundle  $\pi: \sP \to M$ with the structure group $G$ is a vector subbundle $\sH \subset T\sP$ which is invariant under the right $G$-action on $\sP$ and which  splits the exact sequence $$ 0 \longrightarrow {\rm Ker}({\rm d} \pi) \longrightarrow T \sP \longrightarrow \pi^* TM \longrightarrow 0. $$  
\end{definition}
 
 An important feature of a $G$-structure $\sP \subset \F M$ is that it carries a distinguished $V$-valued 1-form.
 
 \begin{definition}\label{d.solder}
 For a $G$-structure $\pi: \sP \to M$, its {\em soldering form} is a $V$-valued 1-form $\theta$ on $\sP$ defined as follows. For each $\varepsilon \in \sP$ and $ x = \pi(\varepsilon)$, the homomorphism $\theta_{\varepsilon}: T_{\varepsilon} \sP \to V$   is the composition
 $$ T_{\varepsilon} \sP \ \stackrel{ {\rm d}_{\varepsilon} \pi}{\longrightarrow} \ T_x M \ \stackrel{\varepsilon^{-1}}{\longrightarrow} \ V,$$ where $\varepsilon^{-1}$ is the inverse of the isomorphism $\varepsilon \in \sP_x \subset \F_x M = {\rm Isom}(V, T_x M)$.  \end{definition}

 Using the soldering form, we can define the following key invariant of a principal connection on a $G$-structure.
 
  \begin{definition}\label{d.torsion}
 Let $M$ be a complex manifold and let $\pi: \sP \to M$ be a $G$-structure with a principal connection $\sH \subset T\sP.$ Then for each vector $u \in V$ and each point  $\varepsilon \in \sP,$ there is a unique  vector $u^{\sH_{\varepsilon}} \in \sH_{\varepsilon} \subset  T_{\varepsilon} \sP$  that satisfies $\varepsilon (u) = {\rm d} \pi(u^{\sH_{\varepsilon}}).$ The {\em torsion function} of the connection $\sH$ is the $\Hom(\wedge^2 V, V)$-valued holomorphic function $\tau^{\sH}$ on $\sP$, whose value $\tau^{\sH}_{\varepsilon}$ at a point $\varepsilon \in \sP$ is  defined by $$\tau^{\sH}_{\varepsilon} (u, v) := 2  \ {\rm d} \theta (u^{\sH_{\varepsilon}}, v^{\sH_{\varepsilon}})$$ for any  $u, v \in V$ in terms of the soldering form $\theta$ on $\sP$.  We say that a  principal connection $\sH$  is {\em torsion-free} if its torsion function vanishes: $\tau^{\sH} =0$.   \end{definition}

  \section{G-structure associated to a $Z$-cone structure}\label{s.GZ}
  How do we use the tools from Section \ref{s.G} to study $Z$-cone structures in Section \ref{s.cone}?
For a projective submanifold $Z \subset \BP V$, let $\h{Z} \subset V$ be its affine cone and let $$G^Z \:= \{ g \in {\rm GL}(V) \mid g(\h{Z}) = \h{Z}\}$$ be its linear automorphism group. 
 Given a $Z$-cone structure $\sC \subset \BP TM$, set for each $x \in M$, 
 $$\sP_x \ := \{ \varepsilon \in \F_x M = {\rm Isom}(V, T_x M)  \mid \varepsilon(\h{Z}) = \h{\sC}_x \}.$$ Then $\sP := \cup_{x \in M} \sP_x \ \subset \ \F M$ is a $G^Z$-structure on $M$ for the group $G^Z \subset {\rm GL}(V).$ Conversely, given a $G^Z$-structure $\sP \subset \F M$ on $M$, set $$ \h{\sC}_x \ := \ \{ \varepsilon(\h{Z}) \subset T_x M \mid \varepsilon \in \sP_x \subset {\rm Isom}(V, T_x M)\}.$$ Then $\sC := \cup_{x \in M} \sC_x \ \subset \ \BP TM$ is a $Z$-cone structure on $M$. In other words, there is a natural one-to-one correspondence:
$$ \mbox{$Z$-cone structures}  \  \Leftrightarrow \ \mbox{  $G^Z$-structures.} $$
Obviously, the  local equivalence of $Z$-cone structures corresponds to the local equivalence of $G^Z$-structures. 
Hence locally flat (resp.  locally symmetric) $Z$-cone structures correspond to locally flat (resp. locally symmetric) $G$-structures. Thus problems like Question \ref{q.Zcone} (hence Question \ref{q.Z}) can be studied using the tools introduced in Section \ref{s.G}. In this section, we recall some general results in this approach.

 \medskip
 Let us say that a $Z$-cone structure is {\em torsion-free} if the associated $G^Z$-structure admits a torsion-free connection in the sense of Definition \ref{d.torsion}. 
If a $Z$-cone structure is locally symmetric, it is torsion-free (see \cite[Chapter XI, Theorem 3.2]{KN}). Sometimes the converse is true. We have the following examples.

\begin{example}\label{ex.Id} Suppose that the Lie algebra $\fg^Z \subset \fgl(V)$ of $G^Z \subset {\rm GL}(V)$ is just $\C {\rm Id}_V$ (note that $\fg^Z$ always contains $\C {\rm Id}_V$). In fact, this is the case for a general choice of a projective submanifold $Z \subset \BP V$. It is not hard to check (for example, it is implicitly proved in \cite[Theorem 3.4]{Hw10})  that, when $\fg^Z = \C {\rm Id}_V$,  a $Z$-cone structure is locally flat if and only if it is torsion-free.  \end{example}
  
 \begin{example}\label{ex.MS}  
Suppose that  $G^Z \subset {\rm GL}(V)$ is an irreducible representation and $Z \subset \BP V$ is the unique closed orbit of the $G^Z$-action on $\BP V$.  We call such a projective submanifold $Z \subset \BP V$ a {\em highest weight variety}. In this case, \cite[Theorem 6.12]{MS} implies that   a torsion-free $G^Z$-structure  is locally symmetric, unless $Z \subset \BP V$ is one of the following.
 \begin{itemize} \item[(1)]  The highest weight variety $Z \subset \BP V$ is the VMRT at a  point of  an irreducible Hermitian symmetric space. Depending on the type of the Hermitian symmetric space, it is one of the following. \begin{itemize} \item[(I)]  The Segre embedding $\BP^{a-1} \times \BP^{b-1} \subset \BP^{ab-1}$.
 \item[(II)] The Pl\"ucker embedding ${\rm Gr}(2; W) \subset \BP \wedge^2 W$ of  the Grassmannian  of 2-dimensional subspaces in a complex vector space $W$.
     \item[(III)] The second Veronese embedding  $v_2(\BP W) \subset \BP \Sym^2 W$ for a vector space $W$.
     \item[(IV)] The smooth quadric hypersurface $\Q^{n-2} \subset \BP^{n-1}.$
\item[(V)] The spinor embedding $\BS_5 \subset \BP^{15}$ of the spinor variety $\BS_5$ of
 isotropic subspaces of dimension 5 in an orthogonal vector space of dimension 10.
\item[(VI)] The 16-dimensional highest weight variety  $Z \subset \BP^{26}$ of the basic representation of the exceptional Lie group $E_6$. \end{itemize}
\item[(2)] Subadjoint varieties:
\begin{itemize} \item[($BD$)] The Segre product  $\BP^1 \times \Q^{m-2} \subset \BP^{2m-1}.$ \item[($G_2$)] The twisted cubic curve $v_3(\BP^1) \subset \BP^3$. \item[($F_4$)]     The Pl\"ucker embedding of the Lagrangian Grassmannian of  isotropic subspaces of dimension 3 in  a symplectic vector space of dimension 6.
\item[($E_6$)]   The Pl\"ucker embedding ${\rm Gr}(3; \C^6) \subset \BP (\wedge^3 \C^6)$  of the Grassmannian  of 3-dimensional subspaces in $\C^6$.
\item[($E_7$)] The spinor embedding $\BS_6 \subset \BP^{31}$ of the spinor variety $\BS_6$ of
 isotropic subspaces of dimension 6 in an orthogonal vector space of dimension 12.
\item[($E_8$)] The 27-dimensional highest weight variety $Z \subset\BP^{55}$ of the basic representation of the exceptional Lie group $E_7$.
\end{itemize}
     \end{itemize} 
      Furthermore, locally symmetric $G^Z$-structures for highest weight varieties $Z \subset \BP V$ are classified up to local equivalences by E. Cartan (see \cite[Chapter XI]{KN} and \cite[Section 26.5]{Ti}). \end{example}
 
 From these examples, we can see that, in many cases of $Z \subset \BP V$,  the torsion functions of principal connections play an essential role in the study of Questions \ref{q.Z} and \ref{q.Zcone}.  
 An important algebraic notion in the study of torsion functions is the following. 
 
 \begin{definition}\label{d.Spencer}
 Let $\fg \subset \fgl(V)= {\rm End}(V)$ be a Lie subalgebra. The {\em Spencer homomorphism} $\p: \Hom (V, \fg) \to \Hom(\wedge^2 V, V)$ is defined by $$\p h(u,v) := h(u) \cdot v - h(v) \cdot u$$ for any $u,v \in V$ and $h \in \Hom(V, \fg)$, where $h(u) \cdot v$ denotes the action of $h(u) \in \fg \subset {\rm End}(V)$ on $v$. \end{definition}
 
 If the Spencer homomorphism is injective for the Lie algebra $\fg \subset \fgl(V)$ of the group $G \subset {\rm GL}(V)$,  then a $G$-structure can have at most one torsion-free principal connection (see \cite[Theorem 2.7]{Sc} or \cite[Proposition 2.8]{HL24c}). This is the case for the Lie algebra $\fg$ of $G^Z$ for most $Z \subset \BP V$, as the following result from \cite[Theorem 7.5]{FH} shows. 
 
 \begin{theorem}\label{t.FH} 
 Let $Z \subset\BP V$ be a projective submanifold. Then  the Spencer homomorphism is injective for the Lie algebra $\fg \subset \fgl(V)$ of $G^Z$, unless
  $Z \subset \BP V$ is projectively equivalent to one in the following list.
\begin{itemize} \item[(a)] The variety in the class (1) of Example \ref{ex.MS}. 
\item[(b)] Let $W$ be a vector spaces of dimension  at least $2.$ Denote by $Q$ a vector space and, by abuse of notation, also a trivial vector bundle with the fiber equal to $Q$ on $\BP W$.  Let $L$ be the tautological line bundle on $\BP W$.  Let $Z$ be the projective bundle $\BP ((Q \otimes L) \oplus L^{\otimes 2})$ over $\BP W$, which is embedded in $\BP V$ by the complete linear system $$H^0(\BP W, (Q \otimes L^*) \oplus (L^*)^{\otimes 2}) = (Q \otimes W)^* \oplus \Sym^2 W^* =: V^*.$$
\item[(c)] A nonsingular linear section with codimension at most two of the Pl\"ucker embedding of the  6-dimensional Grassmannian $\Gr(2; \C^5) \subset \BP^9$.
\item[(d)] A $\BP^4$-general linear section with codimension at most three of the 10-dimensional spinor variety $\mathbb{S}_{5} \subset \BP^{15}$ in Example \ref{ex.MS} (V). 
\item[(e)] Some biregular projections of (a) and (b), which are completely described in  \cite[Section 4]{FH12}.
\end{itemize}  \end{theorem}
 
 In (d), a $\BP^4$-general linear section means that a general element among linear sections containing a linear subspace $\BP^4 \subset \mathbb{S}_{5}$. 
 
The following result from \cite[Theorem 3.9]{HL24c} is useful in the study of Questions \ref{q.Z} and \ref{q.Zcone}.

\begin{theorem}\label{t.connection}
Let $Z \subset \BP V$ be a projective submanifold. For
 a point  $0 \neq v \in \h{Z}$, denote by $T_v \h{Z}$  the affine tangent space of $\h{Z}$ at $v$. Define the subspace $\Xi_Z \subset \Hom(\wedge^2 V, V)$ by
 $$\Xi_Z := \{ \sigma \in \Hom(\wedge^2 V, V) \mid \sigma(v, T_v \h{Z}) \subset T_v \h{Z} \mbox{ for any } 0 \neq v \in \h{Z}\}.$$ 
Let $\fg \subset \fgl(V)$ be the Lie algebra of $G^Z \subset {\rm GL}(V)$. 
Assume that $Z$ satisfies the following three conditions.  \begin{itemize}
\item[(i)] $H^0(Z, \sO(1)) = V^*$;
\item[(ii)] the natural homomorphism $\fg \otimes V^* \to H^0(Z, TZ \otimes \sO(1))$ induced by the restriction homomorphism $\fg \to H^0(Z, TZ)$ is surjective; and
    \item[(iii)] $\Xi_Z$ is contained in the image $\p(\Hom(V, \fg)) \subset \Hom(\wedge^2 V, V)$ of the Spencer homomorphism.  
     \end{itemize}
    Then a  $Z$-cone structure equipped with a characteristic conic connection is torsion-free.  \end{theorem}
 
 Theorem \ref{t.connection} gives a reasonable answer to Question \ref{q.Zcone} in the setting of Examples \ref{ex.Id} and \ref{ex.MS}.  This can be used to answer Question \ref{q.Z} in many cases.

 \section{Results on Question \ref{q.Z}}\label{s.Z}   
 
 In this section, we survey  results on Question \ref{q.Z} for some specific projective varieties $Z \subset \BP V$.
 
\subsection{When $Z \subset \BP^{n-1}$ is a linear subspace}\label{s.linear}

The first case when Question \ref{q.Z} was studied was when  $Z = \BP^{n-1}$, which is the VMRT of projective space $X = \BP^n$. In this case, an affirmative answer is given by the following stronger result of \cite[Theorem 0.2]{CMSB}.

\begin{theorem}\label{t.CMSB}
Let $X$ be a  projective manifold of dimension $n$. If there is a family of minimal rational curves on $X$
such that the VMRT $\sC_x$ at some point $x \in X$ is equal to $\BP T_x X$, then $X$ is projective space $\BP^n$ and minimal rational curves are  lines on $\BP^n$.\end{theorem} 

The argument in \cite{CMSB} is purely algebraic-geometric. This was extended to the case when $Z \subset \BP V$ is a linear subspace in \cite[Theorem 1.1]{Ar} as follows.  

\begin{theorem}\label{t.Araujo}   
In Question \ref{q.Z}, if $Z \subset \BP V$ is a linear subspace, then a general member $C \subset X$ of $\sK$ is of locally flat type.
\end{theorem} 

 In this case, the VMRT $\sC \subset \BP TX$ determines a distribution on a Zariski-open subset of $X$. The proof of Theorem \ref{t.Araujo} in   \cite{Ar} uses mostly algebraic geometry, but also an elementary result from differential geometry: Frobenius theorem applied to this distribution.

\subsection{When $Z \subset \BP V$ is a complete intersection}\label{s.CI}

The most typical example of  $Z \subset \BP V$ with $\fg^Z = \C {\rm Id}_V$ in Example \ref{ex.Id} is a complete intersection of degree $\geq 3$.  This case is studied in \cite[Theorem 1.7]{FH18}, where the following result is obtained. 

\begin{theorem}\label{t.FH18}
In Question \ref{q.Zcone}, assume that $Z \subset \BP V$ is a complete intersection of codimension $c$ and of multi-degree $[m_1, \ldots, m_c]$ satisfying one of the following.
\begin{itemize} 
\item[(I)] $Z$ is a curve of multi-degree different from $[3], [4], [2,2], [2,3]$ or $[2,2,2]$.
\item[(II)] $Z$ is covered by lines and of multi-degree different from $[2], [3]$ or $[2,2]$.
\item[(III)] $Z$ is contained in a nonsingular hypersurface of degree $d \geq 3$ and its multi-degree satisfies $d= m_1 < d+2 \leq m_2 \leq \cdots  \leq m_c.$ \end{itemize}
    Then a $Z$-cone structure with a characteristic conic connection is locally flat. Consequently, the curve $C \subset X$ in Question \ref{q.Z} is of locally flat type for  $Z \subset \BP V$. \end{theorem}
    
This is proved by applying Theorem \ref{t.connection} (actually, its simpler version \cite[Corollary 2.27]{FH18} for $\fg = \C {\rm Id}_V$).  The linear normality condition (i) of Theorem \ref{t.connection}  is always satisfied by a complete intersection. The condition (ii)  can be checked  by standard methods of algebraic geometry. Let us explain how to check the condition (iii). 

From Theorem \ref{t.FH} and  $\fg = \C {\rm Id}_V$, the image of the Spencer homomorphism is isomorphic to 
$$V^* \cong \{ \sigma \in \Hom(\wedge^2 V, V) \mid \sigma(u, v) \subset \C u + \C v \mbox{ for any } u,v \in V\},$$
which is contained in  $$\Xi_Z := \{ \sigma \in \Hom(\wedge^2 V, V) \mid \sigma(v, T_v \h{Z}) \subset T_v \h{Z} \mbox{ for any } v \in \h{Z}\}.$$ Thus, to check the condition (iii), it suffices to prove $\dim \Xi_Z = \dim V$.

To check $\dim \Xi_Z = \dim V,$ we use a homomorphism $$\zeta: \Xi_Z   \longrightarrow  H^0(Z, TZ \otimes T^*Z \otimes \sO(1))$$ defined as follows.
Note that the fiber of $TZ \otimes T^*Z \otimes \sO(1)$ at a point $[v] \in Z$ corresponding to $0 \neq v \in \h{Z}$ is naturally isomorphic to $$\Hom( v  \otimes (T_v \h{Z}/\C v), T_v \h{Z}/\C v).$$
For an element $\sigma \in \Xi_Z$, the value of $\zeta(\sigma) \in H^0(Z, TZ \otimes T^*Z \otimes \sO(1))$ at $[v] \in Z$, as an element of $ \Hom( v  \otimes (T_v \h{Z}/\C v), T_v \h{Z}/\C v),$ is defined to be  $$\zeta(\sigma)_{[v]} (v \otimes (w \mod \C v)) = \sigma(v,w) \mod \C v $$ for any $w \in T_v \h{Z}.$  It is easy to check that $\zeta$ is injective for a nondegenerate complete intersection $Z \subset \BP V$, unless $Z$ is a plane conic. 
In the cases (I) and (II) of Theorem \ref{t.FH18}, one checks $\dim \Xi_Z = \dim V$  by showing that $\dim H^0(Z, TZ \otimes T^*Z \otimes \sO(1)) = \dim V$. For the case (III), one shows that $\dim H^0(Y,TY \otimes T^*Y \otimes \sO(1)) = \dim V$ for the hypersurface $Y$ of degree $d$ containing $Z$ and then proves that 
 the image  $$\zeta(\Xi_Z) \ \subset \ H^0(Z,TZ \otimes T^*Z \otimes \sO(1))$$ lies in the restriction of  $H^0(Y,TY \otimes T^*Y \otimes \sO(1))$ to $Z \subset Y$. 

The excluded cases  of multi-degrees  in Theorem \ref{t.FH18} are due to some technical difficulties in the computations explained above.   

\begin{question} Are the excluded cases in Theorem \ref{t.FH18} necessary? Namely, 
do examples which are not  locally flat exist in  the excluded cases? Do such examples exist in the setting of Question \ref{q.Z}? \end{question}

For example, in the setting of Question \ref{q.Z},  when $Z \subset \BP V$ is a hypersurface of degree $2$ or $3$, which are excluded  in Theorem \ref{t.FH18} (II), it is known that a general member $C \subset X$ of $\sK$ is of locally flat type.  

\subsection{When $Z \subset \BP V$ is an adjoint variety}\label{s.adjoint}

When $Z \subset \BP V$ has the automorphism group  $G^Z$ of dimension bigger than 1, it is often difficult to check the conditions of Theorem \ref{t.connection}, especially (iii).  One case where the conditions have been checked successfully is the adjoint variety of a simple Lie algebra.

Let  $\fg$ be a simple Lie algebra. Choose a complex Lie group $G$ of the Lie algebra $\fg$ such  that the adjoint representation of $G$ is faithful, giving an injection $G \subset {\rm GL}(\fg)$. Its highest weight variety $Z^{\fg} \subset \BP \fg$ is called the {\em adjoint variety} of $\fg$. Regarding $G$ as a complex manifold and $\fg$ as the tangent space $T_e G$ at the identity element $e \in G$, the left translation of $Z^{\fg} \subset \BP T_e G$ determines a $Z^{\fg}$-cone structure $\sC^{\fg} \subset \BP TG$. The subvariety $\sC^{\fg} \subset \BP TG$ is preserved under both left and right translations and the inversion of $G$. Thus this $Z^{\fg}$-structure is locally symmetric. One can show that it is not locally flat, except when $\fg$ is of Dynkin diagram type $C$ (symplectic Lie algebras). 

In \cite{DP}, De Concini and Procesi have constructed a projective manifold $X^G$, called the {\em wonderful compactification} of $G$ such that \begin{itemize}
\item[(a)]  there is a $G\times G$-action on $X^G$ with an open orbit $O^G\subset X^{G}$;
\item[(b)] there is a  biregular morphism from $G$ to $O^G$ which is equivariant with respect to the $G\times G$-action on $G$ by the  left and right translations; and 
     \item[(c)]  the boundary $X^{G} \setminus O^G$ has nice geometric properties described in \cite[page 1-2]{DP}. \end{itemize}
Brion and Fu showed in \cite[Theorem 1.1]{BF}  that unless $\fg$ is of Dynkin diagram type $A$, there is a unique family $\sK$ of minimal rational curves on $X^G$ such that its total VMRT corresponds to $\sC^{\fg} \subset \BP TG$ under the equivariant biregular morphism $O^G \cong G$ in (b). Thus we see that in Question \ref{q.Z} with $Z \subset \BP V$ projectively equivalent to the adjoint variety $Z^{\fg} \subset \BP \fg$ with $\fg$ of Dynkin diagram type neither $A$ nor  $C$, there are minimal rational curves which are not of locally flat type. We have the following classification result from \cite[Theorem 1.7]{HL24c}

\begin{theorem}\label{t.adjoint}
Assume that $Z \subset \BP V$ is the adjoint variety of a simple Lie algebra $\fg$ of Dynkin diagram type neither  $A$ nor  $C$. Then in Question \ref{q.Zcone}, a $Z$-cone structure with a characteristic conic connection is locally symmetric. In Question \ref{q.Z},  a general member $C \subset X$ of $\sK$ is  either of locally flat type or germ-equivalent to a general minimal rational curve on the wonderful compactification $X^G$.  \end{theorem}
      
This is proved by applying Theorem \ref{t.connection}. The conditions (i) and (ii) of Theorem \ref{t.connection} is easy to check for the adjoint varieties. In fact, they hold even when $\fg$ is of type $A$ or $C$. To check the condition (iii), the contact structure of the adjoint variety is used. Each adjoint variety $Z^{\fg}$ carries a natural distribution $D^{\fg} \subset T Z^{\fg}$  
of corank 1, which is a contact structure. The proof is by  examining  the homomorphism $\zeta: \Xi_Z  \to H^0(Z, TZ \otimes T^*Z \otimes \sO(1))$ explained in \ref{s.CI} with respect to this contact structure $D^{\fg} \subset T Z^{\fg}$. See \cite[Sections 5 and 6]{HL24c} for details. Once these conditions are checked, Theorem \ref{t.connection} guarantees that the $Z$-cone structure is torsion-free. Then  \cite[Theorem 6.12]{MS} cited in Example \ref{ex.MS} says that the underlying $G^Z$-structure is locally symmetric, implying the first statement of Theorem \ref{t.adjoint}. For the second statement of Theorem \ref{t.adjoint}, one needs to determine  possible curvature tensors on such locally symmetric structures, which is done in \cite[Section 7]{HL24c}.  

For $\fg$  of Dynkin diagram type $C$, the adjoint variety is projectively equivalent to the Veronese embedding  in Example \ref{ex.MS} (III). In this case, a general minimal rational curve is of locally flat type, as discussed in \ref{s.IHSS} below.

For $\fg$  of Dynkin diagram type $A$, the set of germ equivalence classes of minimal rational curves is infinite. In fact, \cite[Theorem 1.2]{BFM} shows that hyperplane sections of ${\rm Gr}(k; \C^{2k}), k \geq 3,$ give infinitely many distinct Fano manifolds of Picard number 1 whose VMRT at a general point is the adjoint variety of type $A$. Hence their minimal rational curves have distinct germ-equivalence classes by Theorem \ref{t.CF}. 
 
The simple group $G$ as a $G\times G$-homogeneous space is an example of semisimple symmetric spaces. The construction of \cite{DP} gives  wonderful compactifications of  semisimple symmetric spaces.  Brion, Kim and Perrin have generalized \cite{BF} to wonderful compactifications of semisimple symmetric spaces in \cite[Table 1]{BKP}. It raises the following question. 

\begin{question} For which cases in the list of \cite[Table 1]{BKP}, does an analog of Theorem \ref{t.adjoint} hold? \end{question}

\subsection{When $Z \subset \BP^{n-1}$ is the VMRT of an irreducible Hermitian symmetric space}\label{s.IHSS}

For the exceptional cases listed in Example \ref{ex.MS} and Theorem \ref{t.FH}, we need some special approaches to study Question \ref{q.Z}. Ngaiming Mok has done a pioneering work in this direction,
proving the following result in \cite[Section 2]{Mo08}. Mok has not stated it explicitly in \cite{Mo08}, but it is the essential point in the proof of his Main Theorem. 

\begin{theorem}\label{t.Mok}
In Question \ref{q.Z}, if $Z \subset \BP V$ is one of the class (1) in Example \ref{ex.MS}, then a general member $C \subset X$ of $\sK$ is of locally flat type.
\end{theorem}

Let us sketch Mok's argument. In retrospect, it consists of the following two key ingredients. 

\begin{itemize} \item {\bf (Extension)} The geometric structure can be extended along rational curves.
\item {\bf (Canonical connection)} The geometric structure induces a canonical connection on a canonical principal bundle such that the vanishing of the curvature tensors of the connection implies the local flatness of the geometric structure. \end{itemize}

In the setting of Theorem \ref{t.Mok},  `the geometric structure' refers to the $G^Z$-structure on  a Zariski open subset $M \subset X$ given by the total VMRT. The ingredient 
{\bf (Extension)} says that the open subset $M \subset X$ contains general members of $\sK$. 
The proof of this assertion consists of two components of different nature. The first is  the following general result in \cite[Proposition 2.2]{Mo08}.

\begin{proposition}\label{p.Mok}
Let $\sC \subset \BP TX$ be the VMRT of a family $\sK$ of minimal rational curves on  a uniruled projective manifold $X$. Let $C \subset X$ be a smooth  unbendable rational curve belonging to $\sK$ and let $C^{\sharp} \subset \sC$ be the lifted curve consisting of tangents to $C$. Then the isomorphism types of the  second fundamental forms of the family of projective subvarieties $\{ \sC_x \subset \BP T_x X \mid x \in C\}$ along $C^{\sharp}$ are constant (independent of $x \in C$). \end{proposition}

Note that Proposition \ref{p.Mok} is completely general, with no restriction on the VMRT $\sC_x$ at a general point. It is a rather simple consequence of the isomorphism type (\ref{e.N}) of the normal bundle of an unbendable rational curve. 

  The second component in the proof of {\bf (Extension)} is the following lemma.

\begin{lemma}\label{l.FF} 
 Let $\Delta = \{ t \in \C \mid |t| < 1\}$ be the unit disc. Let $\sZ \subset \Delta \times \BP V$ be an irreducible closed analytic subset such that the fiber $$Z_t := \sZ \cap (\{t\} \times \BP V)   \ \subset \ \BP V $$ is projectively equivalent to  $Z \subset \BP V$ for all $t \in \Delta, t \neq 0$, where  $Z \subset \BP V$ is in the class (1) in Example \ref{ex.MS}.  Assume that for a general point $z$ in an irreducible component $Z'_0$ of $Z_0$, the second  fundamental form of $Z'_0 \subset \BP V$ at $z$ is isomorphic to the second fundamental form of a general point of $Z \subset \BP V$. Then $Z'_0 \subset \BP V$ is projectively equivalent to $Z \subset \BP V$ and $Z_0 = Z'_0$. \end{lemma}

In \cite[Proposition 2.3]{Mo08}, Lemma \ref{l.FF} is checked by extending the  G-structure on $Z_t$, which is defined by the second fundamental forms, to $Z'_0$ and then applying the complex-analytic result of Hartogs extension theorem. It is also possible to give a more differential-geometric argument, using \cite[Theorem 1]{HY}.  

By combining Proposition \ref{p.Mok} and Lemma \ref{l.FF}, we see that the VMRT $\sC_x \subset \BP T_x X$ at {\em every} point $x \in C$ of a general member $C \subset X$ of $\sK$ is projectively equivalent to $Z \subset \BP V$, hence $C$ is contained in the open subset $M \subset X$. This establishes {\bf (Extension)}. 

In the setting of Theorem \ref{t.Mok}, the second ingredient {\bf (Canonical connection)} is classical. It was established in \cite{Oc}, or more generally, in \cite[Theorem 4.15]{Ta79}.  Mok's argument  uses this pre-established theory from differential geometry.

When we have both ({\bf Extension}) and ({\bf Canonical connection}), we can  conclude that the $G^Z$-structure is locally flat. In fact, one can easily show that   curvature tensors of holomorphic connections vanish on $M$ by considering their restrictions to general members of $\sK$. 

\medskip
Mok's argument can be generalized to the case when $Z \subset \BP V$ is VMRT of a rational homogeneous space $G/P$, where $P$ is a maximal parabolic subgroup of a simple Lie group $G$ associated with a long simple root. This was done in \cite[Section 3]{Mo08} when $G/P$ is of contact type and in \cite{HH08} for the remaining cases. However, these results do not say that the minimal rational curves are of locally flat type. In fact,  the corresponding projective variety $Z \subset \BP V$ in this case  is linearly degenerate, unless $G/P$ is  an irreducible Hermitian symmetric space. In particular, the VMRT determines a nontrivial distribution on a neighborhood of the minimal rational curve, which is non-integrable on $G/P$.  This means that the corresponding $G^Z$-structure is not locally flat. The results in \cite[Section 3]{Mo08}  and \cite{HH08} say that if, in the setting of Question \ref{q.Z} applied to the linearly degenerate $Z \subset \BP V$ of $G/P$, the successive brackets of the distribution determined by the VMRT span the tangent bundle at general points, then the minimal rational curve is germ-equivalent to a minimal rational curve on $G/P$. 

 There is a considerable difficulty in extending this argument  to the case when $Z \subset \BP V$ is VMRT of a rational homogeneous space associated with a short simple root, one example of which is the symplectic Grassmannian, to be discussed in \ref{s.symplectic} below. 

\medskip
 Let us mention that  an alternate proof of Theorem \ref{t.Mok} is given in \cite{HN}, which does not use the argument {\bf (Extension)} and uses more refined theory of {\bf (Canonical connection)}, the  relative version of parabolic geometry. In this approach, the role of unbendable rational curves is to give some additional local condition on the $Z$-cone structure. This method may be extendible to some  other $Z$-cone structures.

\subsection{When $Z \subset \BP V$ is the VMRT of a (odd) symplectic Grassmannian}\label{s.symplectic}

Let $\Sigma$ be a vector space and $\omega$ be a skew-symmetric 2-form on $\Sigma$ of maximal rank. If $\dim \Sigma$ is even (resp. odd), then $\omega$ is called a {\em symplectic form} (resp. an {\em odd-symplectic form}) on $\Sigma$. Let $X$ be ${\rm Gr}_{\omega}(k, \Sigma),$ the set of all $k$-dimensional $\omega$-isotropic subspaces of $\Sigma$ for $k \geq 2$, called a {\em symplectic Grassmannian} (resp. an {\em odd symplectic Grassmannian}) when $\dim \Sigma$ is even (resp. odd).  The automorphism group of $X$ has a dense open orbit $O \subset X$ ($O = X$ when $\dim \Sigma$ is even). There is a unique family $\sK$ of minimal rational curves on $X$ and the VMRT $\sC_x \subset \BP T_x X$  at a point $x =[W] \in O$, corresponding to a $k$-dimensional isotropic subspace $W \subset \Sigma$, is projectively equivalent to the variety in Theorem \ref{t.FH} (b) for a suitable choice of $Q$. 

It turns out that the a general minimal rational curve on $X$ is {\em not} of locally flat type.  This can be seen as follows. The subspace of $V^*$ corresponding to $\Sym^2 W^*$ in the complete linear system in Theorem \ref{t.FH} (b) determines a distinguished subspace in $T_x X$ at each point $x \in O$, namely, a distribution $D \subset TO$. If a general member $C \subset X$ of $\sK$ is of locally flat type, the associated $G^Z$-structure is locally flat, which implies that $D$ is an integrable distribution. But it is well-known that $D$ is not integrable. 

Thus for $Z \subset \BP V$ of Theorem \ref{t.FH} (b), there are at least two distinct germ-equivalence classes in Question \ref{q.Z}. In fact, we may choose $\omega$ of smaller rank and still consider 
${\rm Gr}_{\omega}(k, \Sigma)$ for suitable values of $k$, called a {\em pre-symplectic Grassmannian}. Then ${\rm Gr}_{\omega}(k, \Sigma)$ has singularity, but we can set $X$ to be a desingularization of ${\rm Gr}_{\omega}(k, \Sigma)$ and  obtain a family of minimal rational curves on $X$ whose VMRT at a general point is projectively equivalent to the same $Z \subset \BP V$. It turns out that if   $\omega$ and $\w{\omega}$ are skew-symmetric forms of different ranks, then the germ-equivalence classes of general minimal rational curves on ${\rm Gr}_{\omega}(k, \Sigma)$ and ${\rm Gr}_{\w{\omega}}(k, \Sigma)$ are different. So we obtain several different germ-equivalence classes in Question \ref{q.Z},  depending on the rank of $\omega$. Their classification is obtained in the  following result from \cite[Theorem 1.5]{HL21}. 

\begin{theorem}\label{t.HL21}
In Question \ref{q.Z}, if $Z \subset \BP V$ is as Theorem \ref{t.FH} (b), then a general member $C \subset X$ is germ-equivalent to a general minimal rational curve on a desingularization of the  pre-symplectic Grassmannian ${\rm Gr}_{\omega}(k, \Sigma)$ for some skew-symmetric form $\omega$ on a vector space $\Sigma$. \end{theorem}

The proof of Theorem \ref{t.HL21} takes a similar approach as that of Theorem \ref{t.Mok}, roughly  consisting of the two ingredients {\bf (Extension)} and {\bf (Canonical connection)} in Subsection \ref{s.IHSS}. 
In fact, it is relatively easy to check an analogue of Lemma \ref{l.FF} for $Z \subset \BP V$ in Theorem \ref{t.FH} (b), using an argument similar to Mok's. Thus the ingredient {\bf (Extension)} can be established in the setting of Theorem \ref{t.HL21}, just as in Theorem \ref{t.Mok}.

The main difficulty lies in the ingredient {\bf (Canonical connection)}. There is no pre-existing differential-geometric theory of $G^Z$-structures in this case and one has to build such a theory. Moreover, the above examples of minimal rational curves on pre-symplectic Grassmannians suggest  that it is not a single structure theory that has to be constructed, but a number of different structure theories,  depending on the properties of the distribution $D$ given by the subspace of $V^*$ corresponding to $\Sym^2 W^*$. A method of constructing canonical connections for such geometric structures subordinate to distributions has been suggested by  \cite{Mo93}, generalizing \cite{Ta79}. However, there is a serious difficulty in applying this method in the current situation.  
 The crucial difference is that while the automorphism group of $Z$ in 
Theorem \ref{t.FH} (a) (hence \cite{Oc} or \cite{Ta79}) is reductive, the automorphism group of $Z$ in Theorem \ref{t.FH} (b) is not reductive. 

This difficulty was resolved in \cite{HL21}, by fusing the two ingredients {\bf (Extension)} and {\bf (Canonical connection)}. A rough sketch of the argument is as follows. First one tries to adapt the prolongation procedure of \cite{Ta79} (more precisely, its simplified version in \cite{CS}) to obtain a canonical principal bundle on $M$ with a canonical connection associated with the geometric structure subordinate to a distribution. This is done in an induction argument and each induction step usually encounters some obstruction.  The key idea is to regard the obstruction in each step  as a holomorphic tensor field on $M$. Then using minimal rational curves lying on $M$ from {\bf (Extension)}, one can show that this obstruction tensor field vanishes on $M$. Thus one can proceed to the next level of the induction. Finally, one obtains the desired canonical connection and its curvature tensors, and prove its vanishing by the restriction to minimal rational curves. 

In short, we may describe the difference between the proof of Theorem \ref{t.Mok} and that of Theorem \ref{t.HL21} as follows. In the proof of Theorem \ref{t.Mok}, {\bf (Extension)} and {\bf (Canonical connection)} are used independently. In particular, the differential geometry in {\bf (Canonical connection)} is a local theory of differential-geometric structures. In the proof of Theorem \ref{t.HL21}, {\bf (Canonical connection)} depends on {\bf (Extension)}, consequently, the differential geometry in {\bf (Canonical connection)} is not a local differential geometry, depending on the global geometry of the neighborhoods of general minimal rational curves from {\bf (Extension)}. 

  The approach of \cite{HL21} is applicable to some other cases of $Z \subset \BP V$. For example, in \cite{HL24} and \cite{HK25}, it is used to answer Question \ref{q.Z} for VMRT of smooth horospherical varieties of Picard number 1, classified by \cite{Pa}. We mention that \cite{HK25} has used a differential-geometric result from \cite{HM24}, which gives a more systematic version of  the prolongation argument in \cite{HL21}.

We expect  that the approach of \cite{HL21} is applicable to the remaining cases (c) and (d) of Theorem \ref{t.FH}, maybe even (e). But the actual implementation in each case seems to require a new idea. A hyperplane section of $\BS_{5}$ in Theorem \ref{t.FH} (d) is the VMRT of a rational homogeneous space with Lie algebra of type $F_4$. In this case, Question \ref{q.Z} is currently being investigated by  Q. Li, D. Timashev and myself.

\subsection{When $Z \subset \BP V$ is a subadjoint variety}\label{s.subadjoint}

We have mentioned that the adjoint variety $Z^{\fg} \subset \BP \fg$ in Subsection \ref{s.adjoint} has a contact structure $D^{\fg} \subset T Z^{\fg}$. When $\fg$ is of Dynkin diagram type neither $A$ nor $C$,  lines lying on $Z^{\fg}$ form a family $\sK$ of minimal rational curves on $Z^{\fg}$ and the VMRT of $\sK$ at a point $x \in Z^{\fg}$ gives a  projective submanifold
$\sC_x  \subset \BP D_x$, which is a highest weight variety. This is called the {\em subadjoint variety} of $\fg$. The class (2) in Example \ref{ex.MS} is the list of subadjoint varieties, labeled by the Dynkin diagram type of $\fg$. 

The answer to Question \ref{q.Z} in this case is the following from \cite[Theorem 1.5]{HL22}. 

\begin{theorem}\label{t.HL22}
In Question \ref{q.Z}, if $Z \subset \BP V$ is one of the class (2) in Example \ref{ex.MS}, excepting the case $(G_2)$ of the twisted cubic curve,  then a general member $C \subset X$ of $\sK$ is of locally flat type. \end{theorem}

We can sketch the argument in terms of the two ingredients {\bf (Extension)} and {\bf (Canonical connection)} in Subsection \ref{s.IHSS}. As before, we have a $G^Z$-structure on a Zariski-open subset $M \subset X$.
Since the Lie algebra $\fg$ of the group $G^Z$ is reductive in this case and moreover the Spencer homomorphism is injective, we can use the standard differential geometric theory to construct a canonical connection on the $G^Z$-structure. Thus the ingredient {\bf (Canonical connection)} looks easy. The difficulty lies in {\bf (Extension)}, because,  unlike in the setting of Theorem \ref{t.Mok} and Theorem \ref{t.HL21}, we can not hope that $M$ contains rational curves. In fact, for $X$ in Proposition \ref{p.blowup} applied to $Z \subset \BP^{n-1}$ of Theorem \ref{t.HL22}, we can check that the corresponding $M \subset X$ contains no rational curves. Thus neither the argument in  \ref{s.IHSS} nor the one in \ref{s.symplectic} could be applied directly here.

The approach in \cite{HL22} is to work on the parameter space  $\sK$ of minimal rational curves, instead of $X$. Because of some special projective-geometric properties of the class (2) in Example \ref{ex.MS}, there is a smooth Zariski-open subset $U \subset \sK$ that carries a natural geometric structure subordinate to a certain distribution. Furthermore, the open subset $U$ contains the $\rho$-images of general fibers of $\mu$ in (\ref{eq.mu}).  In particular, the open set $U$ is covered by submanifolds biregular to $Z$, hence by rational curves.  Thus {\bf (Extension)} is automatic for this geometric structure on $U$. Then the issue is to establish {\bf (Canonical connection)} for this geometric structure. It turns out that the argument of \cite{HL21} can be applied here to give a canonical connection. In other words,  we can use {\bf (Extension)} to construct the canonical connection on $U$, as in Subsection \ref{s.symplectic}. This way, we have the local flatness of the geometric structure on $U \subset \sK$, from which we can deduce the local flatness of the $G^Z$-structure on $M \subset X$ via the universal family (\ref{eq.mu}). 

The reason that the twisted cubic curve is excluded in Theorem \ref{t.HL22} is that the argument of \cite{HL21} cannot be applied to the associated geometric structure on $U$, because the underlying Lie algebra does not satisfy a certain finiteness condition. In fact, we suspect that when $Z \subset \BP V$   is the twisted cubic curve, there may be examples of  minimal rational curves in Question \ref{q.Z}, which are not locally symmetric.

\section{Formal principle for minimal rational curves}\label{s.FP}

In this final section, we discuss a rather different aspect of Question \ref{q.equiv}.
An essential difficulty in Question \ref{q.equiv} is that the existence of a biholomorphic map $\varphi: U \to \w{U}$ is a transcendental question. There is no effective computational method to check it. 
A classical approach of reducing it to an effectively computable problem is to consider  finite-order neighborhoods. Let $(C/X)_k$ (resp. $(\w{C}/\w{X})_k$) be the $k$-th order neighborhood of $C$ in $X$ (resp. $\w{C}$ in $\w{X}$) defined by the $(k+1)$-th power of the ideal of $C \subset X$ (resp. $\w{C} \subset \w{X}$). A necessary condition of the germ-equivalence is the existence of an isomorphism between the nonreduced schemes $(C/X)_k$ and $(\w{C}/\w{X})_k$ for each $k>0$. This is, at least theoretically, a computable algebraic problem for each $k$. 
In fact, the germ-equivalence implies that one can choose the  isomorphisms between $(C/X)_k$ and $(\w{C}/\w{X})_k$  for all $k$ such that they are compatible with each other,  yielding a formal isomorphism between the formal neighborhoods $(C/X)_{\infty} = \lim_k (C/X)_k$ and $(\w{C}/\w{X})_{\infty} = \lim_k (\w{C}/\w{X})_k$. 
A natural question is whether this necessary condition is sufficient for the germ-equivalence. This turns out to be a highly nontrivial question in complex analysis. Let us recall the standard terminology. 

\begin{definition}\label{d.FP} For a compact complex submanifold $A$ in a complex manifold $X$, denote by $(A/X)_{\infty}$ the formal neighborhood of $A$ in $X$. We say that $A \subset X$
 satisfies {\em the formal principle},
    if for any  compact complex submanifold $\widetilde{A}$ in a complex manifold $\widetilde{X}$ admitting
      a formal isomorphism $\psi: (A/X)_{\infty} \stackrel{\cong}{\to} (\widetilde{A}/\widetilde{X})_{\infty}$ between the formal neighborhoods, the submanifolds $A \subset X$ and $\w{A} \subset \w{X}$ are germ-equivalent, namely,  there exists a biholomorphic map $\Psi: U \stackrel{\cong}{\to} \w{U}$ between some  neighborhoods $U \subset X$ of $A$ and $\widetilde{U} \subset \widetilde{X}$ of $\widetilde{A}$ such that $\Psi(A) = \w{A}$.
 \end{definition}

In Definition \ref{d.FP},  the biholomorphic map $\Psi$ may not be an extension of  the formal isomorphism $\psi$. If we require it to be an extension, we have the following  stronger version.

\begin{definition}\label{d.convergence}
A compact complex submanifold $A$ in a complex manifold $X$ satisfies {\em the formal principle with convergence},  if for  any  compact complex submanifold $\widetilde{A}$ in a complex manifold $\widetilde{X}$ admitting
      a formal isomorphism $\psi: (A/X)_{\infty} \stackrel{\cong}{\to} (\widetilde{A}/\widetilde{X})_{\infty}$, 
      there exists a biholomorphic map $\Psi: U \stackrel{\cong}{\to} \w{U}$ between some  neighborhoods $U \subset X$ of $A$ and $\widetilde{U} \subset \widetilde{X}$ of $\widetilde{A}$ such that
      $\psi= \Psi|_{(A/X)_{\infty}},$ giving the following commutative diagram.
  $$ \begin{array}{ccccccc}   X & \supset & U & \stackrel{\Psi}{\longrightarrow} & \widetilde{U} & \subset  & X \\ & & \cup & & \cup & & \\ & & (A/X)_{\infty} & \stackrel{\psi}{\longrightarrow} & (\widetilde{A}/\widetilde{X})_{\infty}. & & \end{array}$$
  In other words, any formal isomorphism $\psi$ is convergent.
\end{definition}

One can find in 
\cite[Section VII.4]{GPR} a good survey of classical results on the formal principle. In the setting of Question \ref{q.equiv}, the following is a special case of Hirschowitz's conjecture from \cite[Section 1]{Hi}.

\begin{conjecture}\label{c.Hi}
A smooth unbendable rational curve satisfies the formal principle.
\end{conjecture}

A partial answer is obtained in \cite[Theorem 1.5]{Hw19}. In the context of Question \ref{q.equiv},  we can state it as follows.

\begin{theorem}\label{t.Hw19}
A general smooth member of any family of minimal rational curves on a uniruled projective manifold satisfies the formal principle.
\end{theorem}

Before discussing the proof, let us ask  whether we can replace the formal principle in the statement of  Theorem \ref{t.Hw19} by the formal principle with convergence.  This is certainly not always possible. Consider a smooth projective curve $R$ and let $X = R \times \BP^1$. Then the $\BP^1$-factor gives a family $\sK = R$ of minimal rational curves on $X$. For any $y \in R,$ let  $C := \{y\} \times \BP^1 \subset X$ be the corresponding minimal rational curve. In a holomorphic coordinate $z$ in a neighborhood of $y$ in $R$,  a formal power series of the form $$f(z) = z + a_2 z^2 + a_3 z^3 +  \cdots$$ gives a formal isomorphism from $(C/X)_{\infty}$ to itself.
This formal isomorphism diverges if the formal power series $f(z)$ is divergent. So any member of $\sK$ violates the formal principle with convergence. This shows that we cannot strengthen Theorem \ref{t.Hw19} to a statement of the formal principle with convergence. 

 In the above example, the VMRT $\sC_x \subset \BP T_x X$ at any $x \in X$ is linearly degenerate in $\BP T_x X$. It turns out that the linear degeneracy is what is behind the failure of the formal principle with convergence. The following theorem is a consequence of a more general result \cite[Theorem 1.2, Proposition 9.3]{HH24}. Some special cases were treated in \cite[Theorem 1.13]{Hw24}.

\begin{theorem}\label{t.HH24}
Assume that the VMRT $\sC_x \subset \BP T_x X$ of a family $\sK$ of minimal rational curves on a uniruled projective manifold $X$ is linearly nondegenerate at a general point $x \in X$. Then a general smooth member of  $\sK$ satisfies the formal principle with convergence.
\end{theorem}

 A concrete example is the following from \cite[Corollary 1.9]{HH24}.

\begin{corollary}\label{c.hypersurface}
A general line on  a smooth hypersurface of degree less than $n$ in $\BP^{n+1}, n \geq 4,$  satisfies the formal principle with convergence, while  a general line on a smooth hypersurface of degree $n$  does not satisfy the formal principle with convergence.
\end{corollary}

The proofs of Theorems \ref{t.Hw19} and \ref{t.HH24} use differential-geometric methods. But it is not through the G-structure given by VMRT on an open subset of $X$. Instead, one should consider a G-structure on the universal family ${\rm Univ}_{\sK}$ in (\ref{eq.mu}).  Let us choose a smooth dense open subset $\sU \subset {\rm Univ}_{\sK}$ such that  the restriction of the universal family morphisms to $\sU$,
$$ \sK \ \stackrel{\rho}{\longleftarrow} \ \sU \ \stackrel{\mu}{\longrightarrow} \ X$$ are two smooth fibrations with ${\rm Ker}({\rm d}_{\alpha} \rho) \cap {\rm Ker}({\rm d}_{\alpha} \mu) =0$ for any $\alpha \in \sU.$
Let $V$ be a vector space of dimension equal to $\dim \sU$ and let $d$ be the fiber dimension of $\mu$.  Fix two subspaces $V_1, V_2 \subset V$ satisfying $\dim V_1 = 1, \dim V_2 = d$  and $V_1 \cap V_2 =0$. Define 
$$G := \{ g \in {\rm GL}(V) \mid g(V_1) = V_1, g(V_2) = V_2\}.$$
Then  we have a $G$-structure $\sP \subset \F \sU$ on $\sU$ with the fiber at $\alpha \in \sU$ given by 
$$\sP_{\alpha} = \{ \varepsilon \in {\rm Isom}(V, T_{\alpha} \sU) \mid \varepsilon(V_1) = {\rm Ker}({\rm d}_{\alpha} \rho), \varepsilon(V_2) = {\rm Ker}({\rm d}_{\alpha} \mu) \}.$$
If a smooth rational curve $\w{C}\subset \w{X}$ is formally equivalent to $C \subset X$,
we have a similar universal family 
$$\w{\sK} \ \stackrel{\w{\rho}}{\longleftarrow} \ {\rm Univ}_{\w{\sK}} \ \stackrel{\w{\mu}}{\longrightarrow} \ \w{X}$$ and a $G$-structure $\w{\sP} \subset \F \w{\sU}$ on a Zariski-open subset $\w{\sU} \subset {\rm Univ}_{\w{\sK}}$. 
It turns out that when $C \subset X$ and $\w{C} \subset X$ are general members of $\sK$ and $\w{\sK}$, respectively, they are germ-equivalent if and only if the $G$-structure $\sP$  at $\alpha \in \sU$ over some point $x \in C$ is locally equivalent to the $G$-structure $\w{\sP}$ at  the corresponding point in $\w{\sU}$, in the sense of Definition \ref{d.Gstructure} (iii).  Moreover, the formal principle and the formal principle with convergence in this setting can be translated as follows.

\begin{question}\label{q.lift}
Assume that  the two $G$-structures $\sP$ on $\sU$ and $\w{\sP}$ on $\w{\sU}$ are formally equivalent at  $\alpha \in \sU$ and $\w{\alpha} \in \w{\sU}$, in  the sense that there is a formal isomorphism $\phi: (\alpha/\sU)_{\infty} \to (\w{\alpha}/\w{\sU})_{\infty}$ that sends $\sP|_{(\alpha/\sU)_{\infty}}$ to $\w{\sP}|_{(\w{\alpha}/\w{\sU})_{\infty}}$. 
\begin{itemize} \item[(i)] Are the $G$-structures $\sP$ at $\alpha$ and $\w{\sP}$ at   $\w{\alpha}$ locally equivalent?
\item[(ii)] Is the formal isomorphism $\phi$ convergent? \end{itemize} \end{question}

It is not difficult to see that an affirmative answer to (i) (resp. (ii)) would imply  the formal principle (resp. the formal principle with convergence) for the corresponding rational curve $C \subset X$.

Theorem \ref{t.Hw19} was proved in \cite{Hw19} by showing that (i) has an affirmative answer if the point $\alpha$ is sufficiently general on $\sU$. The proof uses Cartan's method of equivalence. Roughly speaking, Cartan had proposed a general procedure to check the local equivalence of $G$-structures, which had been rigorously developed by several mathematicians, in particular, \cite{Mo83}. It gives an affirmative answer to (i) at the  point where the system of differential equations describing the solutions of the equivalence problem becomes involutive,
which holds for a sufficiently general point. 

Theorem \ref{t.HH24} was proved in \cite{HH24} by showing that (ii) has an affirmative answer if the point $\alpha$ is sufficiently general on $\sU$ and the $G$-structure $\sP$ satisfies some 
additional conditions (like the one coming from the linear nondegeneracy of VMRT). The proper differential geometric tool in this case is a generalization of Tanaka's prolongation procedure in \cite{Ta70}, which is a much more elaborate version of Cartan's method of equivalence. It is closely related to the ingredient {\bf (Canonical connection)} in Subsection \ref{s.IHSS}. 
The main point is that under  suitable additional conditions, the $G$-structure $\sP$ canonically determines a fiber space  $p: \sP^{\sharp} \to \sU$ with a canonical affine connection $\nabla$ on $\sP^{\sharp}$. Similarly, there is a canonical fiber space $\w{p}: \w{\sP}^{\sharp} \to \w{\sU}$ with a canonical affine connection $\w{\nabla}$ on $\w{\sP}^{\sharp}$. As the fiber spaces are canonically associated to the $G$-structures, the formal isomorphism $\phi$ in Question \ref{q.lift} can be lifted to a formal isomorphism $\phi^{\sharp}$ satisfying the following commutative diagram  at  suitable points $\alpha^{\sharp} \in \sP^{\sharp}$ and  $\w{\alpha}^{\sharp} \in \w{\sP}^{\sharp}$ over $\alpha \in \sU$ and $\w{\alpha} \in \w{\sU}$, respectively: 
$$ \begin{array}{ccc}  (\alpha^{\sharp}/\sP^{\sharp})_{\infty} & \stackrel{\phi^{\sharp}}{\longrightarrow} & (\w{\alpha}^{\sharp}/\w{\sP}^{\sharp})_{\infty} \\
p \downarrow & & \downarrow \w{p} \\ (\alpha/\sU)_{\infty} & \stackrel{\phi}{\longrightarrow} &  (\w{\alpha}/\w{\sU})_{\infty}. \end{array} $$
 Moreover, since the affine connections are canonical, the formal isomorphism $\phi^{\sharp}$ sends $\nabla$ to $\w{\nabla}$. This implies that $\phi^{\sharp}$ is convergent by the following lemma. 

\begin{lemma}\label{l.KN}
Let $(M, \nabla)$ and $(\widetilde{M}, \widetilde{\nabla})$ be two complex manifolds with holomorphic affine connections.  Suppose there exist points $y \in M$ and $\widetilde{y} \in \widetilde{M}$ with a formal isomorphism $\psi: (y/M)_{\infty} \to (\widetilde{y}/\widetilde{M})_{\infty}$ such that $\psi$ sends the restriction  $\nabla|_{(y/M)_{\infty}}$ to the restriction $\widetilde{\nabla}|_{ (\widetilde{y}/\widetilde{M})_{\infty}}.$ Then $\psi$ is convergent.  \end{lemma}

The lemma is proved in \cite[Chapter VI, Theorem 7.2 and its proof]{KN}. The idea is that 
an affine connection $\nabla$ (resp. $\w{\nabla}$) determines in a canonical way a holomorphic coordinate system, the geodesic normal coordinate system, centered at the point $y \in M$ (resp. $\w{y} \in \w{M}$). Then the pull-back of the geodesic normal coordinates at $\w{y} \in \w{M}$ by the formal isomorphism $\psi$ must be the geodesic normal coordinates at $y \in M$. In other words, the formal isomorphism when expressed in these holomorphic coordinates becomes convergent power series.   

\medskip
Let us finish with two open questions.

\begin{question}\label{q.formal}
In Theorems \ref{t.Hw19}, is it necessary to consider only a general member of $\sK$? More generally, is there an example of  a smooth  rational curve that violates the formal principle? \end{question}

We mention that there is an example of a smooth elliptic curve on a surface that violates the formal principle. But no example of a rational curve violating the formal principle is known.  Similarly, one can ask the following.

\begin{question}\label{q.convergence}
In Theorem \ref{t.HH24}, is it necessary to  consider only a general member of $\sK$? Is there an example with linearly nondegenerate VMRT such that  a smooth  rational curve belonging to $\sK$  violates the formal principle with convergence? For example, is there a line on  a smooth hypersurface of degree less than $n$ in $\BP^{n+1}, n \geq 4,$  that violates the formal principle with convergence? \end{question} 
\bibliographystyle{amsalpha}

     \bigskip

Jun-Muk Hwang (jmhwang@ibs.re.kr)

Center for Complex Geometry,

Institute for Basic Science (IBS),

Daejeon 34126, Republic of Korea
\end{document}